\documentclass{primus}

\usepackage{amsmath,amssymb, xcolor, graphicx, comment, csquotes, enumitem, fancyhdr, tikz, pgfplots, subcaption, hyperref}
\PassOptionsToPackage{hyphens}{url}
\pgfplotsset{compat=1.18}
\usepackage[english]{babel} 
\usepackage[sorting=nyt, url = false, doi = false, isbn =false, backend = biber]{biblatex}
\usepackage[symbol]{footmisc}

\usetikzlibrary{arrows.meta}

\allowdisplaybreaks

\title{ENGAGING ACTIVITIES FOR TEACHING LINEAR ALGEBRA\footnote[2]{The authors contributed equally to the paper.}}

\author{Shintaro Fushida-Hardy\\
Department of Mathematics\\
Stanford University\\
Stanford, CA 94305, USA\\
sfh@stanford.edu
\and
Pranav Nuti\\
Department of Mathematics\\
Stanford University\\
Stanford, CA 94305, USA\\
pranavn@stanford.edu
\and
Megan Selbach-Allen\\
Graduate School of Education\\
Stanford University\\
Stanford, CA 94305, USA\\
mselbach@stanford.edu
}


\newcommand{\AmSLaTeX}{$\cal A$\kern-.1667em\lower.5ex\hbox{$\cal
M$}\kern-.125em $\cal S$-\LaTeX}

\definecolor{CBgreen}{HTML}{00A716}
\definecolor{CBpurple}{HTML}{830088}

\newcommand{\CO}{\color{CBgreen}\textbf{COLL}\color{black}}
\newcommand{\EX}{\color{CBgreen}\textbf{EXPL}\color{black}}
\newcommand{\LI}{\color{CBgreen}\textbf{LINK}\color{black}}
\newcommand{\RE}{\color{CBgreen}\textbf{RELA}\color{black}}
\newcommand{\TE}{\color{CBpurple}\textbf{CURR}\color{black}} 
\newcommand{\TI}{\color{CBpurple}\textbf{TIME}\color{black}}
\newcommand{\BA}{\color{CBpurple}\textbf{BACK}\color{black}}

\begin{document}

\makePtitlepage
\makePtitle

\begin{abstract}
This paper discusses several linear algebra activities designed to help enhance students' skills in collaborating, exploring mathematics, and linking together abstract and visual ways of approaching mathematics. Most of these activities are short, accessible, engaging, and easy to incorporate into any classroom. In addition, we discuss some questions instructors can ask themselves to design novel and engaging activities when constrained to teaching from a particular curriculum. 
\end{abstract}

\listkeywords

\section{INTRODUCTION} 
Succeeding in college mathematics requires developing many skills students may not anticipate they will need. As college math instructors, we desire to help students learn such important skills as how to collaborate, how to
explore and read math by themselves, and how to link together different ways of understanding a concept. We also want our students to develop a positive relationship with math. On the other hand, college instructors in introductory classes often face demands from their institution to cover a certain curriculum in a certain amount of time, and these restrictions can make it hard to teach students the skills we believe they need to succeed in college. To make this task even more challenging, students have widely varying backgrounds as they enter college, even as they take the same few introductory classes.

In this paper we share activities from a linear algebra class we designed for a summer bridge-to-college program for incoming first-year students. We focus on innovative, engaging tasks designed to go beyond the dichotomy between addressing students' needs and dealing with institutional constraints. Our discussion of each activity includes a reflection of how students engaged with it. We expect that the tasks we designed will be useful to instructors teaching introductory linear algebra classes and make all the activities available in this footnote\footnote{All activities designed for the class are available at \url{https://github.com/sfushidahardy/SSEA-Linear-Algebra-22/blob/main/README.md}. The website includes linear algebra and non-linear algebra activities in addition to those described in this paper.}. While the tasks we have created are flexible, we also share our activity design process because we recognize that these tasks may need to be modified to address the constraints of various curricula.

\section{GOALS AND CONSTRAINTS}
\subsection{Goals}
When we began designing the class and activities, we had several goals centered around addressing our students' needs. These goals came out of our understanding of effective teaching, and our pedagogical commitments, which were in turn informed by the math education literature.

While there are other important goals not in the following list, we decided to choose a small number of focus points for our course. Given the high degree of overlap between our course goals, and the four pillars of inquiry-based math education laid out by Laursen and Rassmussen \cite{laursen2019prize}, we would label our approach  inquiry-based. We expect many instructors might already agree with these goals and be familiar with the literature informing them, but we have included citations for anyone who wants to further explore the literature. We wanted students to:

\addtolength{\leftmargini}{0.6cm}
\begin{enumerate}
    \item[\CO] Learn how to \textbf{collaborate} in the context of mathematical problem-solving \cite{treisman1985study}.
    \item[\EX] Learn how to \textbf{explore} and read math by themselves, especially because of the fundamental importance of exploration in problem-solving \cite{freeman2014active, polya2004solve}.
    \item[\LI] Be able to \textbf{link} together different ways of approaching a concept, particularly visual and abstract ways of approaching a concept, but also how the concept appears in the textbook, and how it appears in applications \cite{boaler2016seeing, leong2012presenting}.
    \item[\RE] Develop a positive \textbf{relationship} with mathematics through engaging activities, and develop mathematical confidence \cite{gafoor2015high, hubert2014learners, larkin2016hate}.
\end{enumerate}
\addtolength{\leftmargini}{-0.6cm}

Our choice of topics to cover and activities to design were guided by these goals. While the algebra of solving linear systems is important, we decided to focus instead on concepts such as planes (and their extensions, subspaces) and orthogonality for a couple of reasons:

\paragraph{\EX} Students have some pre-existing familiarity with planes and orthogonality, and these concepts have familiar applications. This makes it easier to prompt students to explore mathematical concepts and applications themselves.

\paragraph{\LI} Understanding these concepts deeply involves bridging the more concrete, visual mathematics many students have experienced in earlier courses with the more abstract descriptions needed to extend ideas to higher dimensions. These notions thus provide excellent material to help students begin linking together different ways of understanding a concept.

\subsection{Constraints}

Like all instructors, we also faced various constraints:

\addtolength{\leftmargini}{0.6cm}
\begin{enumerate}
    \item[\TE] We had to teach from a \textbf{curriculum} over which we did not have full control.
    \item[\TI] We had only so much \textbf{time}, so we had to design activities that could fit into narrow windows.
    \item[\BA] We had to consider that students came with widely varying mathematical \textbf{backgrounds}, and so we had to create activities that could engage students at multiple levels.  
\end{enumerate}
\addtolength{\leftmargini}{-0.6cm}

As instructors we often have freedom in the specifics of how we teach our classes. However, eliminating a topic from the curriculum, or not providing it classroom time might not be an institutionally viable option (especially if the class is a part of a coordinated course sequence), even if we judge that it would be best for our students' mathematical understanding. We also faced this constraint, and had to adapt our lessons to the topics designated by the department for the course.

This constraint was particularly important in our design choices for this class. A primary contribution of this paper is thinking about how we can design interesting activities while being constrained to covering particular topics within a syllabus and taking only a certain amount of time.

\section{ACTIVITY DESIGN}

Our design process started with a couple of basic assumptions, which we believe many instructors will share.
\vspace{-4mm}
\paragraph{\RE} Every activity we designed had to engage the students and help them develop mathematical confidence.
\vspace{-4mm}
\paragraph{\TE} Every activity was based on the curriculum we were constrained to, and the concepts we needed to cover.

The first step in our activity design process was to come up with a rough idea for a task related to the concept we wanted to cover. For some concepts we quickly came up with an initial idea on our own. Other concepts lent themselves to finding an activity online that we could build from. In this case, we would often have to modify the activity to suit the topic we needed to cover. For example, one of the activities we considered was the ‘Carpet Ride’ problem series from the Inquiry Oriented Linear Algebra (IOLA) curriculum \cite{wawro2012inquiry}. For us however, an important topic in the curriculum was the discussion of various ways of describing planes. We thus had to modify the problem series as we discuss later in the paper.

In either case, task development for us was an iterative process. After we had an initial idea in mind, we considered the questions below to repeatedly modify the activity to better fit our goals.

\addtolength{\leftmargini}{0.6cm}
\begin{enumerate}
\item[\CO] Is this activity collaborative? 
\item[\EX] Is there a step in the activity where I'm just telling students something?
\item[\LI] Is there a step of the activity that involves only abstract definitions? Is there some visual component to this task? Is the activity about a real-world application?
\item[\TI] How long do I think this activity will take?
\item[\BA] Does the activity have a simple entry point while simultaneously challenging my students?  Is there a step of the activity that seems like a big leap?
\end{enumerate}
\addtolength{\leftmargini}{-0.6cm}

After every lesson, we discussed how effective activities were and what issues they had. Driving our activity design is our belief that effective instruction is fundamentally student-centered.  We taught this course three consecutive summers, and asked students what aspects of the course they found most helpful.  We revised the overall curriculum and individual activities based on this feedback, which was sought through anonymous surveys and one-on-one interviews. An example of how one activity we discuss in this paper has recently been modified is detailed later in this paper.

In the next section we share a sampling of the activities we developed, and how each activity connects with our goals and the questions listed above. We also discuss how students engaged with each activity.

\section{SAMPLE ACTIVITIES} 

In this section, we share some tasks we designed. Additional tasks are available at \url{https://github.com/sfushidahardy/SSEA-Linear-Algebra-22/blob/main/README.md}. To give a picture of the range of activities we designed, we have selected activities that emphasize different goals and constraints.
\begin{itemize}

\item \emph{Maps and Distances} is designed for students to practice working with magnitude and distance. It also helps them understand the role the origin plays in defining vectors, and how this is related to subtraction.

\item \emph{Perpendicularity and Approximations} is designed to help students understand the importance of perpendicular vectors.

\item \emph{Various Sneaky Mathematicians} is an extension of the IOLA Carpet Ride problem \cite{wawro2012inquiry}, designed to help students discover planes in $\mathbb{R}^3$ and different ways to specify them.

\item \emph{Planes Telephone} is designed to give students practice working with planes and their various descriptions. It also gives students practice with cooperating with others outside of their group.

\end{itemize}

\subsection{Maps and Distances}

\begin{figure}[ht]
    \centering
    \includegraphics[width=0.8\textwidth]{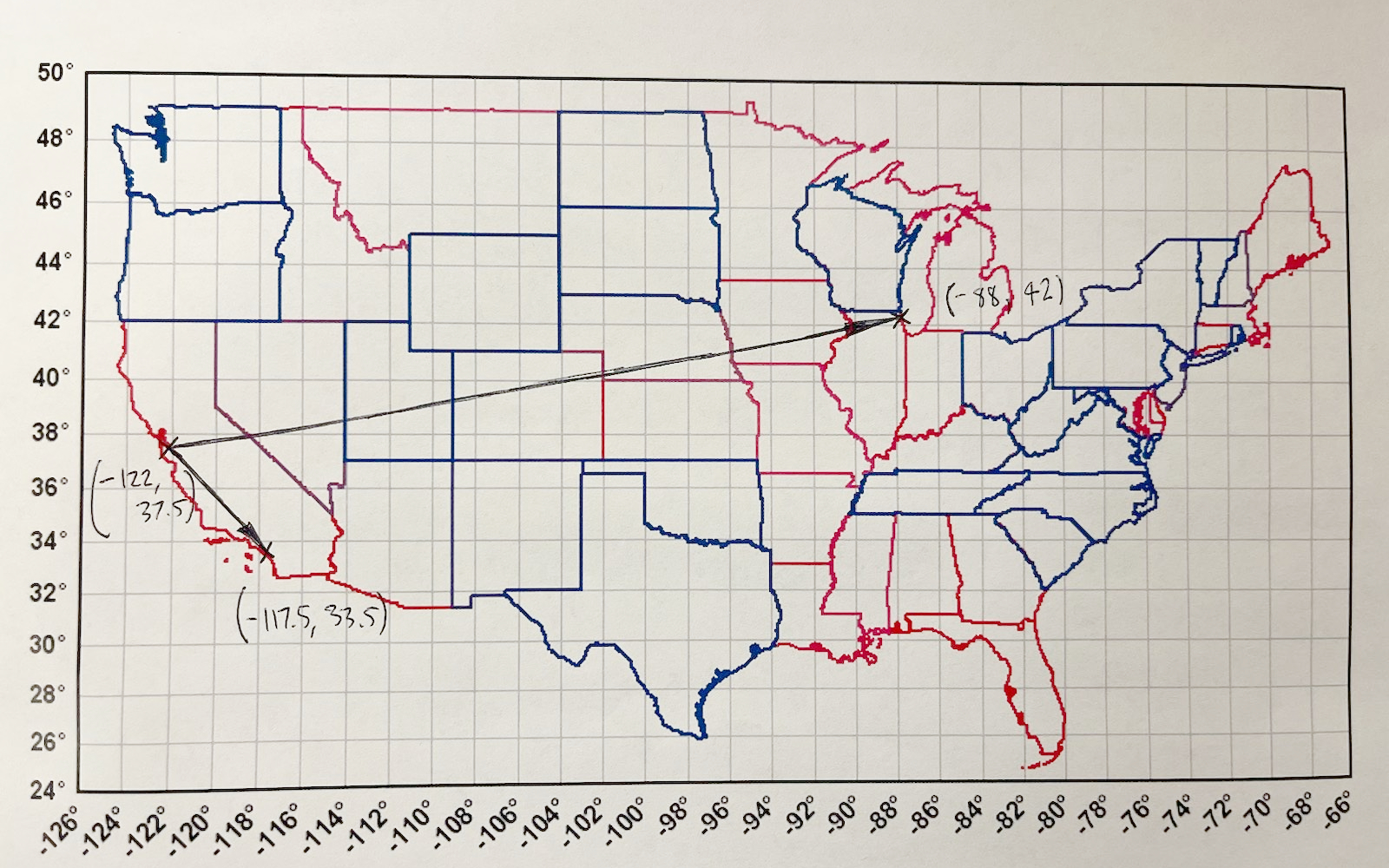}
    \caption{\emph{Maps and Distances} example.}
    \label{fig:map}
\end{figure}

In this activity pairs of students are given maps of California, US, or the world with grid-lines for longitude and latitude. They then proceed in four steps:
\begin{enumerate}
\item Estimate vectors from the university (as the origin) to each of their hometowns. Find the difference of the two vectors and compute its magnitude to approximate the distance between their hometowns. See Figure \ref{fig:map}.
\item Next, set one of their hometowns as the origin and find the vector from this hometown to the other hometown. Use this vector to approximate the distance between hometowns.
\item Compare the two different vector methods of estimating distance between hometowns.  Compare the results with the distance provided by mapping software.
\item Discuss any similarities or differences between the vector distance and the mapping software distance---what could be the cause of any discrepancies?
\end{enumerate}

\subsubsection*{Implementation and reflections}

\paragraph{\CO, \RE} This activity was designed to require collaboration---both students must locate their hometowns and compute the distance between them, which stimulates communication. Moreover, hometowns are typically an important part of a student's identity, especially in the context of a summer bridge program with students who are generally leaving home for the first time. The activity encouraged students to share stories about themselves or their childhoods and connect with each other, enriching their relationships with mathematics.

\paragraph{\LI} Representing vectors on a map is very visual, as is identifying the relationship between the said vectors and the distance between hometowns on the map. On the other hand, computing the distance required students to express the vectors algebraically as elements of $\mathbb{R}^2$ (for example, with components corresponding to longitude and latitude). The students linked together the visual and abstract facets of the problem.

\paragraph{\EX}
Each pair of students received their own unique problem, and thus their own unique answer. Some answers were very accurate while others were not. These features of the activity make it very rich for student exploration.

We did not give strict guidance to students about which units to use in their calculations. This lead to some students counting each square on the map as ``one" while others looked up conversions from longitude and latitude to miles or kilometers. This ambiguity resulted in occasional confusion as students compared answers with each another, but also led to productive discussions of the complexity of these conversions.

In particular, students consistently discovered that some of their calculations were very similar to distances computed by mapping software, while other calculations were much further off. Discussions around this discrepancy highlighted the complexity of converting spherical distances to distances on a flat grid. The discussions helped reveal how different results did not necessarily mean that students had calculated incorrectly, but instead related to issues of scale. Students then discussed whether or not the application of vectors to maps is reasonable, for scales varying between California maps and world maps.

\subsection{Perpendicularity and Approximations} 

\begin{figure}[ht]
  \begin{subfigure}{0.45\textwidth}
    %% Creator: Inkscape 1.2.2 (1:1.2.2+202305151915+b0a8486541), www.inkscape.org
%% Accompanies image file 'perp_and_app_a.pdf' (pdf, eps, ps)
\begingroup%
  \makeatletter%
  \providecommand\color[2][]{%
    \errmessage{(Inkscape) Color is used for the text in Inkscape, but the package 'color.sty' is not loaded}%
    \renewcommand\color[2][]{}%
  }%
  \providecommand\transparent[1]{%
    \errmessage{(Inkscape) Transparency is used (non-zero) for the text in Inkscape, but the package 'transparent.sty' is not loaded}%
    \renewcommand\transparent[1]{}%
  }%
  \providecommand\rotatebox[2]{#2}%
  \newcommand*\fsize{\dimexpr\f@size pt\relax}%
  \newcommand*\lineheight[1]{\fontsize{\fsize}{#1\fsize}\selectfont}%
  \ifx\svgwidth\undefined%
    \setlength{\unitlength}{143.03749275bp}%
    \ifx\svgscale\undefined%
      \relax%
    \else%
      \setlength{\unitlength}{\unitlength * \real{\svgscale}}%
    \fi%
  \else%
    \setlength{\unitlength}{\svgwidth}%
  \fi%
  \global\let\svgwidth\undefined%
  \global\let\svgscale\undefined%
  \makeatother%
  \begin{picture}(1,0.73647752)%
    \lineheight{1}%
    \setlength\tabcolsep{0pt}%
    \put(0,0){\includegraphics[width=\unitlength,page=1]{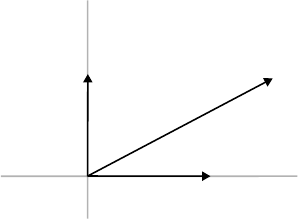}}%
    \put(0.23203872,0.29577131){\makebox(0,0)[lt]{\lineheight{1.25}\smash{\begin{tabular}[t]{l}$u$\end{tabular}}}}%
    \put(0.43648548,0.09541492){\makebox(0,0)[lt]{\lineheight{1.25}\smash{\begin{tabular}[t]{l}$v$\end{tabular}}}}%
    \put(0.57721934,0.34197281){\makebox(0,0)[lt]{\lineheight{1.25}\smash{\begin{tabular}[t]{l}$w$\end{tabular}}}}%
    \put(0,0){\includegraphics[width=\unitlength,page=2]{perp_and_app_a.pdf}}%
  \end{picture}%
\endgroup%

    \caption{Problem 1 in \emph{Perpendicularity and Approximations}.} \label{fig:perp_a}
  \end{subfigure}
  \hspace*{\fill}
  \begin{subfigure}{0.45\textwidth}
    %% Creator: Inkscape 1.2.2 (1:1.2.2+202305151915+b0a8486541), www.inkscape.org
%% Accompanies image file 'perp_and_app_b.pdf' (pdf, eps, ps)
\begingroup%
  \makeatletter%
  \providecommand\color[2][]{%
    \errmessage{(Inkscape) Color is used for the text in Inkscape, but the package 'color.sty' is not loaded}%
    \renewcommand\color[2][]{}%
  }%
  \providecommand\transparent[1]{%
    \errmessage{(Inkscape) Transparency is used (non-zero) for the text in Inkscape, but the package 'transparent.sty' is not loaded}%
    \renewcommand\transparent[1]{}%
  }%
  \providecommand\rotatebox[2]{#2}%
  \newcommand*\fsize{\dimexpr\f@size pt\relax}%
  \newcommand*\lineheight[1]{\fontsize{\fsize}{#1\fsize}\selectfont}%
  \ifx\svgwidth\undefined%
    \setlength{\unitlength}{143.03749275bp}%
    \ifx\svgscale\undefined%
      \relax%
    \else%
      \setlength{\unitlength}{\unitlength * \real{\svgscale}}%
    \fi%
  \else%
    \setlength{\unitlength}{\svgwidth}%
  \fi%
  \global\let\svgwidth\undefined%
  \global\let\svgscale\undefined%
  \makeatother%
  \begin{picture}(1,0.73647752)%
    \lineheight{1}%
    \setlength\tabcolsep{0pt}%
    \put(0,0){\includegraphics[width=\unitlength,page=1]{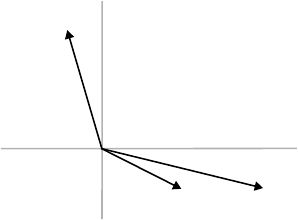}}%
    \put(0.21517861,0.40806359){\makebox(0,0)[lt]{\lineheight{1.25}\smash{\begin{tabular}[t]{l}$w$\end{tabular}}}}%
    \put(0.3933574,0.13703304){\makebox(0,0)[lt]{\lineheight{1.25}\smash{\begin{tabular}[t]{l}$u$\end{tabular}}}}%
    \put(0.63891226,0.17168555){\makebox(0,0)[lt]{\lineheight{1.25}\smash{\begin{tabular}[t]{l}$v$\end{tabular}}}}%
    \put(0,0){\includegraphics[width=\unitlength,page=2]{perp_and_app_b.pdf}}%
  \end{picture}%
\endgroup%

    \caption{Problem 5 in \emph{Perpendicularity and Approximations}.} \label{fig:perp_b}
  \end{subfigure}
\caption{Two problems from \emph{Perpendicularity and Approximations}. Students approximate the vector $w$ as a linear combination of $u$ and $v$.} \label{fig:perp}
\end{figure}

In this activity students (in groups of three or four) are given pictures of vectors $u$, $v$ and $w$, where $\{u,v\}$ span $\mathbb{R}^2$ and range between perpendicular and near-parallel. Students are asked to quickly approximate $w$ as a linear combination of $u$ and $v$ using intuition and visual reasoning. For example, a student might guess in Figure \ref{fig:perp_a} that $w$ is approximately $u + 1.5v$, and in Figure \ref{fig:perp_b} that $w$ is approximately $-4u + v$. After about five minutes, every group records their approximations on a shared table. The class discusses the table, with a focus on the questions ``which pairs of vectors $\{u, v\}$ resulted in similar approximations for $w$ across all groups?" and ``which pairs of vectors had high variation between groups?" Including the discussion, the activity takes about 15 minutes from beginning to end.

\subsubsection*{Implementation and reflections}

\paragraph{\EX} The activity was unique in that its very purpose was to help students uncover a primary application of the concept (orthogonality) by themselves. In the discussion section of the activity, students found that the approximations obtained by the class generally agreed for $u,v$ near-orthogonal, but differed wildly when $u$ and $v$ were near-parallel. See Figure \ref{fig:perp} for examples of $u,v$ near-orthogonal and near-parallel. This led to students discovering that orthogonality is a property that makes decomposition into components more robust. Another observation that students had was that for $u, v$ near-orthogonal, their approximations for $w$ had higher confidence in the sense that the coefficients $a$ and $b$ in the linear combination $au + bv = w$ were given with higher precision. Conversely for $u, v$ near-parallel, these coefficients were typically integers. 

Given more time, the class discussion could delve deeper into \emph{why} orthogonality (or near-orthogonality) of $u$ and $v$ helps to approximate $w$. The instructor may elaborate that orthogonality is also beneficial to computers just as it was to them, and that there are results about expressing vectors as linear combinations later in the course that require orthogonality as well. 

We found during this activity that students sometimes struggled with the task of providing quick, intuitive answers, whether or not $u,v$ were near-orthogonal or near-parallel. This might be because students are used to performing calculations with clear steps when tasked with obtaining numerical answers. We believe that making quick guesses and estimating are extremely useful problem-solving skills that are not stressed enough in many classes. The authors wish to explore this further and create more activities with an emphasis on quantitative intuition.

\paragraph{\TI} The activity only takes about 15 minutes, about half of which is students estimating linear combinations and reporting their results, and the other half is discussion. This activity is well-suited to a preliminary discussion of the topic of orthogonality.

We also found that the activity was a good way to quickly recapitulate the concept of linear combinations (as we moved on to discussing orthogonality), while stressing its geometric nature.

\subsection{Various Sneaky Mathematicians}

This activity is an extension of the Carpet Ride problem from the IOLA project \cite{wawro2012inquiry}. The activity consists of five inquiry-based problems which allow the student to discover planes in $\mathbb{R}^3$ and different ways of describing them. Figures \ref{fig:sneak1}, \ref{fig:sneak2}, and \ref{fig:sneak4} depict problems 1, 2, and 4 of \emph{Various Sneaky Mathematicians} (VSM)---the full activity, including problems 3 and 5, can be found at: \url{https://github.com/sfushidahardy/SSEA-Linear-Algebra-22/blob/main/README.md}. Students are tasked to work on whiteboards or butcher paper, creating posters of their reasoning to share with the rest of the class.

\begin{figure}[ht]
    \centering
    \noindent\fbox{%
    \parbox{\textwidth}{%
    \footnotesize{
        Old Man Gauss is back in $\mathbb{R}^2$ and hidden somewhere along the line $y=7x-23$.  You have a brand new hoverboard and when you start it for the first time you must set the direction of travel (i.e. vector). Once the direction is set it cannot be changed without great difficulty (so you are stuck with the direction you choose).  You also have a dodgy, enchanted portal that can take you to any point in $\mathbb{R}^2$, but you can only be confident it will work once.
\begin{enumerate}
	\item How should you program your hoverboard and where should you take the portal to ensure you can find Old Man Gauss?  
	\item Ask your instructor to learn the exact point $(x,y)$ where Old Man Gauss is hiding. Based on the point and direction you chose in (1), when will you reach Old Man Gauss? 
\end{enumerate}
}
    }%
}
    \caption{Problem 1 of \emph{Various Sneaky Mathematicians}.}
    \label{fig:sneak1}
\end{figure}

\begin{figure}[ht]
    \centering
    \noindent\fbox{%
    \parbox{\textwidth}{%
    \footnotesize{
       Unlike Old Man Gauss, Takakazu Seki has managed to hide in $\mathbb{R}^3$. He is at a point $(x,y,z)$ that satisfies the equation $y=7x-23$.
\begin{enumerate}
	\item Can you reach Seki with the hoverboard and portal you had from Question 1? If your answer is yes, show why this is true.  If your answer is no, what other transport would you need to reach Seki?
\end{enumerate}
}
    }%
}
    \caption{Problem 2 of \emph{Various Sneaky Mathematicians}.}
    \label{fig:sneak2}
\end{figure}

\begin{figure}[ht]
    \centering
    \noindent\fbox{%
    \parbox{\textwidth}{%
    \footnotesize{
       Maryam Mirzakhani is also in $\mathbb{R}^3$ and is hidden at a point satisfying the equation $3x+2y-z=4$.  You have a new hoverboard and magic carpet ready to be programmed (and again you can only program each once).  You also have a single use portal available to jump anywhere in $\mathbb{R}^3$.
\begin{enumerate}
	\item How should you program your hoverboard and magic carpet and how should you use your portal to ensure you can find Mirzakhani?  (Note: remember we are now in $\mathbb{R}^3$ so you will need to specify 3-dimensions for travel points.)
	\item Ask your instructor to learn the exact point $(x,y,z)$ where Mirzakhani is hiding.  Based on the point and direction you chose in (1), when will you reach her?
\end{enumerate}
}
    }%
}
    \caption{Problem 4 of \emph{Various Sneaky Mathematicians}.}
    \label{fig:sneak4}
\end{figure}

\subsubsection*{Implementation and reflections}

\emph{VSM problem 1} (Figure \ref{fig:sneak1}) builds off of the Carpet Ride problem \cite{wawro2012inquiry}. In the Carpet Ride problem, Old Man Gauss is hiding in $\mathbb{R}^2$ and students determine how to find him with two modes of transport (two vectors). In VSM problem 1, students work with a point and a vector, and discover that any point on the line $y = 7x - 23$ in $\mathbb{R}^2$ can be expressed as $p + tv$ where $p$ and $v$ are a fixed point and vector respectively, and $t$ is some ``time" parameter. Indeed, the broad goal of VSM problem 1 is for students to discover that a line in $\mathbb{R}^2$ can be encoded in these different ways.

We saw that students found the openness of the problem challenging, but they grew more comfortable with this openness in later problems. We elaborate on this in the {\EX} section below. 

The broad goal of \emph{VSM problem 2} (Figure \ref{fig:sneak2}) is for students to discover the notion of planes in $\mathbb{R}^3$, and more specifically to convert again from an equational description to a parametric description. Students are provided with the same equation $y = 7x - 23$ as in VSM problem 1, but now must interpret it in $\mathbb{R}^3$. 

In our experience, students successfully recognise that this equation has no $z$ dependence, and ultimately that it describes a \emph{plane} in $\mathbb{R}^3$. They understood that the plane is vertical, and typically thought of it as many copies of the \emph{line} $y=7x-23$ in $\mathbb{R}^2$, but vertically stacked along the $z$-axis. Typical solutions were ``You need an elevator!" or ``You need a rocket!".

\emph{VSM problem 4} (Figure \ref{fig:sneak4}) is another extension of problem 2. This time, the given equation depends on all coordinates $x,y,z$. Students are again (essentially) tasked with converting the equational description to a parametric description.

Our students found this a lot more difficult than VSM problem 2. They seemed to make the connection between $z$-coordinate independence and verticality without much difficulty, but given a general plane they typically found one vector lying in the $xy$-plane and then had difficulty finding a second vector for the parametric form.

\paragraph{\EX} The activity was intentionally designed to avoid the use of words such as \emph{plane}. Students must explore the problem, typically through drawing, using graphing software (which we encouraged), or algebra, to develop an understanding of the relevant spaces, their dimensions, and how they interact. In our experience, students found the activity difficult for two key reasons:
\begin{itemize}
\item They're used to positions being fully unconstrained or fully constrained. (For example, a hypothetical version of problem 1 stating that Old Man Gauss is anywhere in $\mathbb{R}^2$, or stating his exact position.)
\item They're used to questions in mathematics having exactly one answer.
\end{itemize}

In problem 1, many students requested to know where Old Man Gauss was hiding to help with part 1 of the problem. In VSM problem 2 and onwards, our students became a lot more comfortable with partial constraints. Similarly, the many different possible choices for the portal location in problem 1 part 1 added a layer of difficulty for the students, as some had a preconceived idea that only one of the locations could be valid. Again, throughout the activity, students began to recognize that \emph{anything that works is a solution}; there is no one canonical solution.

\paragraph{\RE} Part of the activity involved creating a poster in a group, and later seeing posters by other groups. This gave the students a lot of creative agency, and built their confidence as they worked together to communicate their ideas in playful ways. Creating posters tends to be time consuming, but the freedom was well received. 

The Carpet Ride Problem is contextualized with the names Cramer and Gauss, who both contributed to the development of linear algebra. In VSM we wished to introduce a more diverse group of mathematicians to draw attention away from the perception that mathematics (and linear algebra specifically) was developed in the West. Noteably, Seki was a Japanese mathematician credited with discovering the determinant which is a central concept in linear algebra.  On the other hand, Mirzakhani was a mathematician at the same institution as our course, further humanizing the mathematical experience. Given the context that Cramer and Gauss were mathematicians, some students were naturally curious about these less familiar names, and went out of their way to ask about them or look up their names online. This activity affords the opportunity to use a wide range of names to provoke conversation about the history of mathematics and the humans who have contributed to it.

We also found that the next topic in the course (subspaces, spans, and dimensions) was typically challenging for students. To ease this transition from planes to subspaces, the Various Sneaky Mathematicians activity has been modified in the current version of the course. The analogy between planes and subspaces is now further stressed by considering only planes through the origin, and the ideas of dimension and linear independence are introduced. The handout for this modified activity is provided in the appendix to this paper.

\subsection{Planes Telephone}

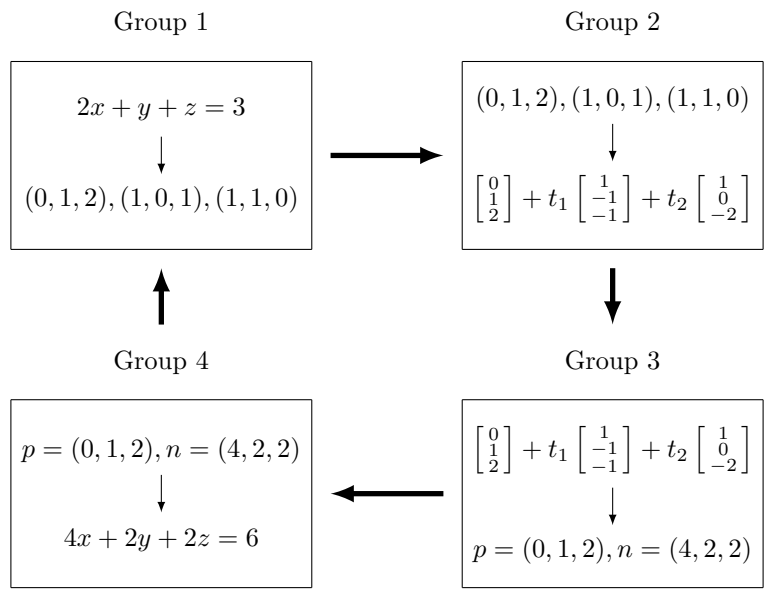
\begin{figure}
    \centering
    \begin{tikzpicture}
		\node [] (0) at (-5, 3.5) {};
		\node [] (1) at (-1, 3.5) {};
		\node [] (2) at (-1, 1) {};
		\node [] (3) at (-5, 1) {};
		\node [] (4) at (1, 3.5) {};
		\node [] (5) at (5, 3.5) {};
		\node [] (6) at (5, 1) {};
		\node [] (7) at (1, 1) {};
		\node [] (8) at (-1, -1) {};
		\node [] (9) at (-5, -1) {};
		\node [] (10) at (-5, -3.5) {};
		\node [] (11) at (-1, -3.5) {};
		\node [] (12) at (1, -1) {};
		\node [] (13) at (1, -3.5) {};
		\node [] (14) at (5, -1) {};
		\node [] (15) at (5, -3.5) {};
		\node [] (16) at (-3, 2.87) {$2x+y+z=3$};
		\node [] (17) at (-3, 4) {Group 1};
		\node [] (18) at (3, 4) {Group 2};
		\node [] (19) at (-3, -0.5) {Group 4};
		\node [] (20) at (3, -0.5) {Group 3};
		\node [] (22) at (3, 3) {$(0, 1, 2) , (1, 0, 1), (1, 1, 0)$};
		\node [] (23) at (3, 1.67) {$\left[\begin{smallmatrix}0\\1\\2\end{smallmatrix}\right] + t_1 \left[\begin{smallmatrix}1\\-1\\-1\end{smallmatrix}\right] + t_2\left[\begin{smallmatrix}1\\0\\-2\end{smallmatrix}\right]$};
		\node [] (24) at (3, -1.67) {$\left[\begin{smallmatrix}0\\1\\2\end{smallmatrix}\right] + t_1 \left[\begin{smallmatrix}1\\-1\\-1\end{smallmatrix}\right] + t_2\left[\begin{smallmatrix}1\\0\\-2\end{smallmatrix}\right]$};
		\node [] (25) at (3, -3) {$p = (0, 1, 2), n = (4, 2, 2)$};
		\node [] (26) at (-3, -1.67) {$p = (0, 1, 2), n = (4, 2, 2)$};
		\node [] (27) at (-3, -2.87) {$4x + 2y + 2z = 6$};
		\node [] (29) at (-3, 1.67) {$(0, 1, 2) , (1, 0, 1), (1, 1, 0)$};
		\node [] (30) at (-0.75, 2.25) {};
		\node [] (31) at (0.75, 2.25) {};
		\node [] (34) at (3, 0.75) {};
		\node [] (35) at (3, 0) {};
		\node [] (38) at (-3, 0) {};
		\node [] (39) at (-3, 0.75) {};
		\node [] (42) at (0.75, -2.25) {};
		\node [] (43) at (-0.75, -2.25) {};
		\node [] (44) at (-3, 2.5) {};
		\node [] (45) at (-3, 2) {};
		\node [] (46) at (3, 2.67) {};
		\node [] (47) at (3, 2.17) {};
		\node [] (48) at (3, -2.17) {};
		\node [] (49) at (3, -2.67) {};
		\node [] (50) at (-3, -2) {};
		\node [] (51) at (-3, -2.5) {};
		\draw (0.center) to (1.center);
		\draw (1.center) to (2.center);
		\draw (2.center) to (3.center);
		\draw (3.center) to (0.center);
		\draw (4.center) to (5.center);
		\draw (5.center) to (6.center);
		\draw (6.center) to (7.center);
		\draw (7.center) to (4.center);
		\draw (12.center) to (14.center);
		\draw (14.center) to (15.center);
		\draw (15.center) to (13.center);
		\draw (13.center) to (12.center);
		\draw (8.center) to (9.center);
		\draw (9.center) to (10.center);
		\draw (10.center) to (11.center);
		\draw (11.center) to (8.center);
		\draw[-latex,line width=2pt] (30.center) to (31.center);
		\draw[-latex,line width=2pt] (34.center) to (35.center);
		\draw[-latex,line width=2pt] (42.center) to (43.center);
		\draw[-latex,line width=2pt] (38.center) to (39.center);
		\draw[-latex] (44.center) to (45.center);
		\draw[-latex] (46.center) to (47.center);
		\draw[-latex] (48.center) to (49.center);
		\draw[-latex] (50.center) to (51.center);
\end{tikzpicture}
    \caption{An example of a plane being converted and communicated in \emph{Planes Telephone}.}
    \label{fig:tele}
\end{figure}

This is another activity relating to descriptions of planes in $\mathbb{R}^3$. This time the class separates into four groups, and plays \emph{telephone}\footnote{Telephone is a classic children's game in which a phrase is whispered from neighbour to neighbour, and laughter ensues when the original and final messages are compared.} with planes---only, each group converts the plane into a new form before communicating it to the next group.
\begin{enumerate}
    \item Each of the four groups is given a \textit{different} plane described in equational form (i.e., an equation in three variables). The group then works together to convert the description to three-point form (i.e., three non-collinear points that lie on the plane), then communicates the result to the group on their left.
    \item When a group receives a plane in three-point form, they convert it to parametric form and communicate the result to their left.
    \item Similarly, parametric form is converted to point-normal form (i.e., a point that lies on the plane, and a vector perpendicular to the plane), and point-normal form back to equational form. When every group has finished with all of the conversions, they receive a plane described in equational form---namely the same plane they started off with! An example of a plane in Planes Telephone is given in Figure \ref{fig:tele}.
\end{enumerate}

\subsubsection*{Implementation and reflections}

\paragraph{\CO, \BA} In this activity the whole class collaborates together on a large scale, but students are also working together on smaller scales. The activity was not competitive, and students with different mathematical backgrounds were all involved in the activity without feeling too much pressure to perform because any mistakes were anonymized through the \emph{telephone} framework. In the traditional game of telephone, the whole point is that the original and final messages differ, which helps reduce pressure.

In our experience, frequently at least one plane was miscommunicated at some point during the activity. This became a further learning opportunity, as the class worked together to determine a step in the conversion where a mistake was introduced. We attempted to value the previously established anonymity by focusing the attention on the \emph{conversion} rather than the students doing the converting. In one of the author's classes, a group discovered that the information they were given couldn't be converted into a valid plane form. Through a discussion, the class learned that this was because an earlier conversion into three-point form actually consisted of collinear points, and thus did not uniquely specify a plane.

Conversely, even when groups completed all of the conversions correctly, their final equation was frequently a multiple of their original equation, rather than being identical.  This discrepancy led to productive discussions about how, even when restricting to a specific form, descriptions of a given plane are not unique.

\paragraph{\TE} This is a topic that is an important part of the linear algebra curriculum at our institution, but students often find the conversions to be involved and conceptually challenging. The activity helped students recognize that the same plane can be represented by different forms. Furthermore, the activity allowed students to practice the skills needed to convert between descriptions of planes, but turned it into a collaborative game which is much more enjoyable for students.

\paragraph{\TI} In our experience the activity takes 20 to 30 minutes, if students have practiced converting between descriptions of planes at least once. This makes it extremely convenient to incorporate into a traditional section or recitation.

\section{CONCLUSION}
To understand student perspectives, we interviewed 10 of the 70 students who took the course in Summer 2021. In analyzing these interviews, we see ways the activities met the goals we felt were important as instructors, and some of the limitations of our work.   

Multiple students reported that they found the activities enjoyable, and one student mentioned how they encouraged creativity. Several students also mentioned that they left the course feeling better about math. One student shared how he uses visuals in a new way because of his experience in the class.  Reviews of student work from the course show multiple examples of student drawings next to algebraic representations.  

Students reported that they still struggled to risk being wrong, which limited some of their mathematical exploration. At the same time, we saw many students making mistakes in class, and we sought to foster an environment where mistakes were encouraged. Indeed, students grew more comfortable making mistakes and working in unfamiliar situations, as described in \emph{Very Sneaky Mathematicians}. 

One goal of the activities was to encourage student collaboration. The overwhelming majority of students reported feeling positively about their experience in group work. Several students mentioned how group work helped them build community, while one reported that group work was intimidating, but not in a bad way. 

However, some students questioned our emphasis on group work, reporting that in their classes during the academic year, they were expected to work alone. They suggested that our course should have more activities that aren't centered around group work. This disconnect points to a need to communicate more to students the value of collaboration, not just for classes in university, but for however they choose to use mathematics afterwards. 

Our interviews revealed that students found our actvities more challenging than their courses during the academic year, which we expected, given the discomfort inquiry-based tasks can cause \cite{deslauriers2019measuring}. Despite this, students also shared how the material covered was helpful in their transition to the academic year, even for the students who were not immediately taking a linear algebra course. 

For most instructors, using more active learning is better for students, and this paper offers ready-to-use activities for teachers of linear algebra. The reflective questions we discussed may spark ideas for designing additional activities. Every instructor faces different constraints, and instructors should feel empowered to design and adapt activities to their curricula and goals.

\section*{ACKNOWLEDGEMENTS}

We would like to express our deepest gratitude to Sophie Libkind, who was an integral member of the syllabus design and teaching team in 2020. Many of her ideas and contributions persist in the present day syllabus of the summer bridge course. We would also like to thank Talia Blum, Julia Costacurta, and Zhihan Li for continuing to improve the curriculum. We would like to offer thanks to Lourdes Andrade and Bethlehem Aynalem, each of whom directed the summer bridge program during our time and trusted and supported us to design the mathematics course. We thank also the many members of the summer bridge program's staff who made the program possible. Finally, we offer thanks to our students over the past three summers for their enthusiasm in engaging with our classes.

\printbibliography[title=REFERENCES]

\section*{APPENDIX}
\label{appendix}
\begin{comment}
\documentclass[11pt, notitlepage, letterpaper]{article}
\usepackage[utf8]{inputenc}
\usepackage{amsfonts, amsmath, amsthm, amssymb, mathtools, xcolor, hyperref, graphicx, enumerate, dsfont, csquotes, paralist, bbm, url, multicol, multirow, wasysym}
\usepackage[english]{babel}
\usepackage[nice]{nicefrac}
\usepackage[margin=1in]{geometry}
\end{comment}
\title{Various Campus Mathematicians}
\author{Talia Blum, Julia Costacurta, Zhihan Li, Pranav Nuti}
\date{SSEA 2023}

\maketitle

\begin{center}
  This activity builds on \href{https://www.tandfonline.com/doi/abs/10.1080/10511970.2012.667516}{IOLA materials on introductory linear algebra.}
\end{center}

You're at $(0,0,0)$ in your dorm room. You have a lot to do. The first two things you have planned are to meet your friend Emmy at Arrillaga for breakfast at $(1,7,4)$, and then you'll head to class in the math building at $(5,20,0)$.

You can take your bike and your skateboard. You're still learning how to ride them, and haven't figured out how to turn without falling. You need to pick a vector direction for each your bike and skateboard, and once you do, you can't change them.

\begin{enumerate}
    \item Which vectors should you set your bike and skateboard to travel along so you can get to both Arrillaga and to your math class?

    \item Given that you have a full schedule of classes for the rest of your day, you decide to check that all of them are actually in places you can travel to on your two vehicles. Calculate all of the locations in 3D space that you can travel to using only your bike and skateboard (set to the directions you chose above). (Hint: it might be helpful to think geometrically.)
\end{enumerate}

\noindent You're all ready to go, but you realize your bike has a flat tire! 

\begin{enumerate}
    \item[3.] Can you travel to both Arrillaga and to class using only your skateboard?

    \item[4.] Describe all of the locations you can travel to on your skateboard. 
\end{enumerate}
    
\noindent Your roommate offers you her scooter, but she warns that it can only travel along the vector $\begin{bmatrix} 0\\3\\4 \end{bmatrix}$. 
\begin{enumerate}
    \item[5.] Can you travel to Arrillaga and to class using the scooter and your skateboard (set to the direction you already chose in part 1)?
\end{enumerate}

\noindent You're \textit{finally} going to head out the door, but you meet your neighbor in the hallway. She offers to fix your flat tire, so you can travel on your bike in the direction you originally intended (from part 1). You'll be late to meet Emmy at Arrillaga if she fixes your tire, so you should only accept her offer if it would increase the number of places you can go on campus. In other words, you should first check if you can reach additional places on campus with all three modes/directions of transit, compared to just having your scooter and skateboard.

\begin{enumerate}
    \item[6.] Should you accept her offer?
    \item[7.] How might this relate to the magic carpet problem? What seems similar or different?
\end{enumerate}

    % You still have a whole list of things to do and places to be for the rest of the day. To make sure you'll be able to make it, what are all the points on campus that you can reach with your bike and skateboard travelling along these vectors?

% \section*{Question 1}
% You're meeting your best friend Emmy for lunch today at the cloud cafe, which is located at the point $(10,20,40)$. You have 3 forms of transit: your bike can travel along the vector $\begin{bmatrix}
%     10 \\ 0 \\0
% \end{bmatrix}$, your skateboard travels along the vector $\begin{bmatrix}
%     0 \\ 2 \\ 4
% \end{bmatrix}$, you can use your rollerblades to travel along $\begin{bmatrix}
%     1 \\ 1 \\ 2
% \end{bmatrix}$.

% \begin{enumerate}
%     \item If you bring your bike, your skateboard, your rollerblades, can you reach the cloud cafe? How?

%     \item Bringing 3 sets of wheels is a lot, so you want to see if you can leave any of your modes of transit at home. Can you get to the cloud cafe using only one or two of your transit options?

% \end{enumerate}

% \noindent After you leave, you remember you have a library book that's due today! You plan to go home to pick up your book, and then hit the library right after seeing Emmy at the cloud cafe. 

% \begin{enumerate}
%     \item Where can the library be located so that you can reach it using your bike, your skateboard, and your rollerblades?
% \end{enumerate}

\section*{BIOGRAPHICAL SKETCHES}

Shintaro Fushida-Hardy grew up in New Zealand and Japan, but recently moved to the United States of America to pursue a PhD in mathematics at Stanford University. Mathematically he is interested in four dimensional topology, and he is also interested in shapes and forms outside of mathematics; he is a sculptor in his free time.

\vspace*{.3 true cm} \noindent Pranav Nuti is a mathematics PhD student at Stanford University, studying sequential decision making. He is interested in mathematics instruction and practices to promote an equitable classroom. He enjoys puzzles, silly improv, and the sea breeze.

\vspace*{.3 true cm} \noindent Megan Selbach-Allen is a PhD candidate in mathematics education at Stanford University. She studies community college mathematics instruction, inquiry-based and active learning pedagogies and alternative assessment practices. Prior to beginning her PhD she taught math at the US Naval Academy. When not working she enjoys weightlifting, walking dogs and flying on the trapeze.

\end{document}